\documentclass{amsproc}[10pt]
\usepackage{amsmath,amscd,amsthm,amssymb,amsxtra,latexsym,epsfig,epic,graphics}
\usepackage[matrix,arrow,curve]{xy}
\usepackage{graphicx}

\def\antiddot{\mathinner{\mkern1mu\raise1pt\vbox{\kern7pt\hbox{.}}\mkern2mu
        \raise4pt\hbox{.}\mkern2mu\raise7pt\hbox{.}\mkern1mu}}

\newcommand{\AAA}{{\mathbb A}}

\newcommand{\CC}{{\mathbb C}}

\newcommand{\FF}{{\mathbb F}}

\newcommand{\PP}{{\mathbb P}}
\newcommand{\QQ}{{\mathbb Q}}
\newcommand{\RR}{{\mathbb R}}

\newcommand{\ZZ}{{\mathbb Z}}

\newcommand{\punkt}{\hspace{-.3ex}\raise.15ex\hbox to1ex{\Huge.}}

\theoremstyle{definition}

\def\CC{{\mathbb C}}
\def\RR{{\mathbb R}}
\def\QQ{{\mathbb Q}}

\makeatletter
\def\Ddots{\mathinner{\mkern1mu\raise\p@
\vbox{\kern7\p@\hbox{.}}\mkern2mu
\raise4\p@\hbox{.}\mkern2mu\raise7\p@\hbox{.}\mkern1mu}}
\makeatother

\makeatletter
\def\Ddots{\mathinner{\mkern1mu\raise\p@
\vbox{\kern7\p@\hbox{.}}\mkern2mu
\raise4\p@\hbox{.}\mkern2mu\raise7\p@\hbox{.}\mkern1mu}}
\makeatother

\date{January 16, 2009}
\title{A family of exceptional Stewart-Gough mechanisms of genus~7}
\author{Florian Gei\ss}
\address{Mathematik und Informatik, Universit\"at des Saarlandes, Campus E2 4, 
D-66123 Saarbr\"ucken, Germany}
\email{fg@math.uni-sb.de}
\author{Frank-Olaf Schreyer} 
\address{Mathematik und Informatik, Universit\"at des Saarlandes, Campus E2 4, 
D-66123 Saarbr\"ucken, Germany}
\email{schreyer@math.uni-sb.de}

\dedicatory{dedicated to Andrew Sommese on the occasion of his sixtieth birthday.}

\begin{document}

\maketitle

\hfill{\it Eppur si muove -- Galileo}

\begin{abstract}
In this paper we construct a family of  exceptional Stewart-Gough Mechanisms, whose motion curves  are algebraic curves of genus 7. Up to translations, rotations and dilatations this family of mechanisms is  $13$-dimensional.

\end{abstract}

\section*{Introduction}

We call a mechanism exceptional if its degree of freedom exceeds the dimension of a general mechanism of the same underlying topological type. In this paper we establish the existence of exceptional Stewart-Gough platforms by a combination of algebraic geometric and number theoretic methods. 

A Stewart-Gough platform consists of two rigid bodies, one of which is grounded, linked to each other by six telescoping legs in ball joints at both ends of each leg (see Figure 1).

\begin{figure}[h]
\begin{center} \includegraphics{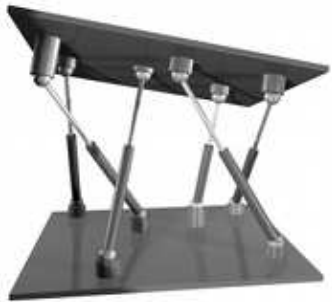} \end{center}
\caption{A general Stewart-Gough platform, Sommese-Wampler \cite{SW}, pp.  107.}
\end{figure}

Usually, for given leg lengths and rigid bodies, there is only a finite number of ways to assemble  the parts.
In robotics,  the ungrounded body is  moved around by changing the leg lengths. 
Since the group of movements $SE(3)$ is six-dimensional and since, in general, when its length is locked, each leg places one constraint on the motion, six legs is the minimum number to rigidly constrain the ungrounded 
body. As the leg lengths are changed, the ungrounded body can be moved 
through a region corresponding to an open part of $SE(3)$.

In this paper we are interested in a different situation: We wish to find a configuration of leg lengths and bodies such that the mechanism moves without changing the leg lengths.  We call such a mechanism exceptional.

Apart from  rather trivial examples, one family of exceptional mechanisms is well-known and thoroughly studied, the Griffis-Duffy platforms \cite{GD}, Husty and Karger \cite{HK1}. The family of Griffis-Duffy platforms forms a unirational variety, Karger \cite{K}. In this paper we ask whether there are families of exceptional Stewart-Gough platforms of different kind. In our main result, we establish  the existence of a family of mechanisms whose motion curve in $SE(3)$ is a smooth algebraic curve of genus 7. In Figure 2, we display the coupler curve, i.e., the curve traced out by a coupler point fixed to the ungrounded link under the motion.

\begin{figure}[h]
\begin{center}
\includegraphics{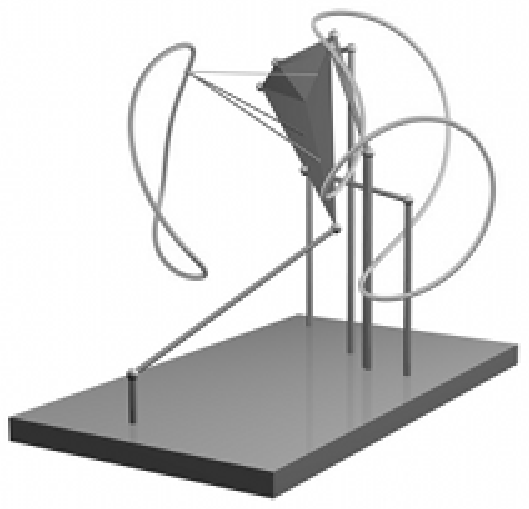}\includegraphics{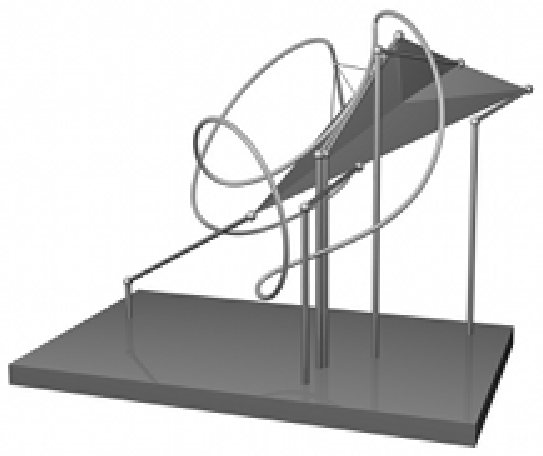}
\end{center}
\caption{A Stewart-Gough platform with a coupler curve of degree 12 and genus $7$.}
\end{figure}


Our methods, we believe, are quite novel in this kind of research. 
We first approach the question by searching for mechanisms defined over a finite field. We then utilize the  algebraic-geometric properties of the discovered mechanisms, to get a construction over  some algebraic number fields $K$, i.e., fields $K$ with $\QQ \subset K \subset \CC$ of degree $[K:\QQ] := \dim_\QQ K < \infty$. Finally, to establish the existince of mechanisms defined over $\RR$, we prove that some of these fields $K$ have  embeddings $K \hookrightarrow \RR$. The final existence result can be verified also by numerical methods as was  pointed out to us by Charles Wampler.

{\bf Acknowledgment.}  We thank Charles W. Wampler and Andrew J. Sommese for introducing us to the problem. We thank Charles for his many  suggestions to improve  the manuscript. We thank Mike Stillman for pointing out a crucial omitted option in an early version of our construction. The second author would like to thank the IMA at the University at Minnesota for  hospitality
and support. The work presented here started during the special program Applications of Algebraic Geometry at the IMA in September 2006.

\section{The algebraic setup}

A movement in the euclidean 3-space is a map
$$\RR^3 \to \RR^3,\, x \mapsto Ax+b$$
with an orthogonal matrix $A=(a_{ij}) \in SO(3)$ and a translation vector $b \in \RR^3$.
We denote the group of movements by $SE(3)$. 

A position of a Stewart-Gough platform is determined by the positions of the six legs. Let $q^{(i)} \in \RR^3$ for $i=1,\ldots,6 $ denote six points on the grounded body, and let $p^{(i)} \in \RR^3$ be the end points on the moving  body. Consider the equations
$$ eq^{(i)}: \qquad \parallel p^{(i)}-q^{(i)} \parallel^2=\parallel Ap^{(i)}+b-q^{(i)} \parallel^2\qquad \hbox{for } i=1,\ldots,6,$$
where $\parallel x \parallel^2 = \left\langle x, x \right\rangle$ denotes the square of the euclidean norm for $ x \in \RR^3$.
Let $$\Phi = V(eq^{(i)} , i = 1, . . . , 6)  \subset \AAA^{36} \times SE(3)$$ be the complete parametrized family of Stewart-Gough mechanisms and 
their assembly positions. Consider the natural projections $\pi_1 :\Phi \to \AAA^{36}$ 
and $\pi_2 : \Phi \to SE(3)$.

Then $C_\alpha = \pi_2 (\pi_1^{-1}(\alpha))$ is the asssembly space of the mechanism specified by $\alpha=(p^{(1)},q^{(1)},\ldots,p^{(6)},q^{(6)}) \in \AAA^{36}$.
For generic $\alpha$, $C_\alpha$ is a finite set of isolated points, 
i.e., $\dim C_\alpha= 0$. We are interested in finding exceptional mechanisms $\alpha$ 
such that $\dim C_\alpha > $0. In particular, when $\dim C_\alpha = 1$, we call the 1-dimensional
part of $C_\alpha$ the motion curve of mechanism $\alpha$. For such $\alpha$ the assembly space
consists of the motion curve and possibly a finite number of isolated points.

The equations $A^T A = E_3$ and $\det A =1$ defining $SE(3)$, and the equation $eq^{(i)}$ make sense as polynomial equations over arbitrary fields $K$. We call $\alpha \in \AAA^{36}(K)$ such that the solution space in 
$C_\alpha \subset SE(3,\overline K)$  is one-dimensional, an exceptional mechanism defined over $K$.
Here $\overline K$  denotes the algebraic closure of $K$. For fixed $\alpha \in \AAA^{36} (K)$ one can determine   the dimension of the solution space  in $SE(3,\overline K)$ with a single Gr\"obner basis computation over $K$. {\sc Macaulay~2} \cite{M2}, which we used, and {\sc Singular} \cite{GPS01} are good computer algebra systems to perform such computations. 

The complete  computer algebra code used is documented online in \cite{documentation}.

\section{Experimental exploration}

Let $M \subset \AAA^{36}$ be a codimension $c$ irreducible  component of the constructible set of exceptional mechanisms.  Heuristically, for the finite field $\FF_q$ with $q$ elements 
the ratio of points on $M$ is approximately
$$\frac{|M(\FF_q)]}{|\AAA^{36}(\FF_q)|}\approx q^{-c}.$$
By the Weil formula, this equation is asymptotically correct for $q \to \infty$ in case of an absolutely irreducible component $M$ defined over $\ZZ$.  

We however plan to make experiments over very small finite fields. 
We implement a script for parallel search for the computer algebra system {\sc Macaulay~2} \cite{M2} running on a cluster with $24$ nodes with $2.4$ GHz each node. 
To compute the dimension of the solution space in $SE(3,\FF_q)$ for one platform takes about $0.3$ seconds. Thus, for very small fields $\mathbb{\Bbbk} = \FF_p$ with $p \in \{3,5,7\}$ we are able to test $p^{N_p}$ platforms, $N_3 = 12$, $N_5 = 9$ and $N_7 = 7$,  each in a maximum time of about $3$ hours. \\ 
We have reasonable hope to find points on $M(\FF_3)$ for components $M \subset \AAA^{36}$ of codimension $\le 12$. To give an idea of such exploration we give the result for a typical test as described above. We list the dimension,  the degree and the total number of the discovered exceptional mechanisms of that degree. Here we use the naive embedding $SE(3)= SO(3) \times \AAA^3 \subset \AAA^{12} \subset \PP^{12}.$
In Table 1
below, the degree reverse to the degree of $C_\alpha$ with respect to the naive embedding.
To get an interpretation for the  total numbers, we also print the hypothetical codimensions, which are only correct, if the family is absolutely irreducible. In case of a family consisting of several components, say $k$, of the same codimension the values $hcodim$ are more likely  approximations of the product $h=\log_p k$ and the codimension $codim$.
\begin{table}[h]\label{expTab}
{\small
$$\begin{array}{rr||rr|rr|rr}
\dim & \deg  & \ \ \  \FF_3 & \ \ \ hcodim & \ \ \  \FF_5 & \ \ \ hcodim &\ \ \  \FF_7 & \ \ \  hcodim \\  \hline \hline
1 & 1  &   2817 & 4.77 & 861 & 4.80  & 60  & 4.90 \\ 
 & 2  &   6223 & 4.05  &1553 & 4.43 & 85  & 4.72 \\
 & 3  &    425 & 6.49  &  20 & 7.14 &  -  & - \\
 & 4  &   2256 & 4.97  & 542 & 5.09 & 29  & 5.27 \\
 & 5  &    383 & 6.59  &  35 & 6.79 &  2  &  \\
 & \bf 6  &   \bf 1083 &  \bf 5.64  &  \bf 580 & \bf 5.05 &  \bf 66  & \bf 4.85 \\
 & 7  &     49 & 8.46  &   2 & & -   & - \\
 & 8  &    236 & 7.53  &   5 &   8.00    & -   & - \\
 & 9  &     12 & 9.74  &   1 &        & -  & - \\
 & 10 &     23 & 9.15  &   - & -       & -  & - \\ 
 & 11 &      2 &   &   1 &        & -  & - \\
 & \bf 12 &     \bf 43 & \bf 8.58  &  \bf  4 & \bf 8.14 & -  & - \\
 & 14 &      5 &  &   - & -       & -  & - \\ 
 & 15 &      1 &        &   - & -       & -  & - \\
 & 16 &      2 &   &   - & -       & -  & - \\
 & 18 &      1 &        &   - & -       & -  & - \\
 & 24 &      3 &        &   - & -       & -  & - \\
 & 28 &      2 &   &   - & -       & -  & - \\ 
 & 32 &      - & -        &   1 &       & -  & - \\
 & 36 &      - & -        &   1 &        & -  & - \\ \hline
2 &  4 &     8 &          &   - & -       & -  & - \\ 
 &  6 &      9 &        &   4 & 8.14& -  & - \\ 
 &  8 &      1 &          &   - & -       & -  & - \\ \hline\hline
& \mathrm{sum} &   13484  &  & 3160 & & 242 
\end{array}$$ 
}\caption{Experimental data}
\end{table}
 
 In the following we investigate the examples of degree $6$ or $12$ further, leaving all other cases to future research. The main reason to choose these cases is that they have  very nice Betti tables, see below.
 To get projective curves, we work in the projective closure $\overline{SE(3)} \subset \PP^{12}$
 of $SE(3) \subset \AAA^{12}$.

We first remove zero-dimensional components of the assembly space, i.e., we compute the homogeneous ideal of the one-dimensional part of the assembly space. This is a standard task in computational algebraic geometry, which for example can be performed using primary decomposition. 
Since this is fairly expensive and more than we want, we actually use an alternative approach built around a probabilistic algorithm, see our documentation \cite{documentation}. The resulting ideals might have  new linear equations, which we can use to eliminate some of the variables. Geometrically, it means that the projective closure of the motion curve $\overline C_\alpha \subset \overline {SE(3)} \subset \PP^{12}$ spans only a subspace $\PP^n$ of smaller dimension.  

The easiest available further information are the Betti numbers of the free resolution of the coordinate rings of the $\overline C_\alpha \subset \PP^n$. These numerical invariants determine the degree and  the arithmetic genus of the curve. Frequently one can deduce the connectedness of the curve, which implies irreducibilty in case of smoothness.  For further information about Betti numbers, we refer to Eisenbud \cite{Eisenbud} and Schreyer \cite{Schreyer}.

The majority of curves of degree $6$ leads to curves $C \subset \PP^5 \subset \PP^{12}$ whose syzygies in $\PP^5$ have Betti table 
 $$
 \begin{matrix}
 1 & -  &  - & - &- \cr
 - & 9 & 16 & 9 &- \cr 
 - & - & - & -& 1 \cr 
 \end{matrix}
 $$
 In particular, the motion curve $C$ is an arithmetically Cohen-Macaulay curve in these cases. The Cohen-Macaulay property of $C$ implies that if $C$ is smooth, then  $C$ is also irreducible.
In summary, $C$ is an elliptic normal curve of degree $6$ in these cases. Unfortunately, none of these curves is the reduction of a curve defined over a real number field; see Gei\ss \ [2008] for an explanation and a unirational description of this family. This shows that in general one cannot deduce the existence of real exceptional mechanisms from the existence of corresponding families over finite fields.
 
 In case of degree $12$ we obtain  curves $C \subset \PP^6$ with Betti table
 $$
 \begin{matrix}
 1 & -  &  - & - &- &-\cr
 - & 10 & 16 & 9 &- &-\cr 
 - & - & 9 & 16& 10 & - \cr
   - & -  &  - & - &- &1\cr
 \end{matrix}
 $$
 Again the motion curve is arithmetically Cohen-Macaulay, hence irreducible if smooth. In this case $C$ is a curve of genus $g=7$ in its canonical embedding. The fact that these curves occur in their canonical embedding is the main reason why we did not work in the  Study quadric $SE(3) \subset Q   \subset \PP^7$ instead of the naive embedding.
 
By Schreyer \cite{Schreyer}, the special shape of the Betti table implies that either $C$ has   a plane model of degree $6$ with 
three (possibly infinitesimally near) double points, or that $C$ is  bi-elliptic. Indeed, the projection onto $\overline{SO(3)} \subset \PP^9$ yields a plane model:
 
 The projection of $C$ to $\overline{SO(3)} \subset \PP^9$ yields a curve $C'$ of degree $12$ which is the intersection of a Veronese surface $V^4$ with a cubic hypersurface $F^3$:
 $$ C' =V^4 \cap F^3 \subset \PP^5 \subset \PP^9$$
 In particular, $C$ has a plane model $C' \subset \PP^2$ of degree $6$. 
 A general plane sextic curve of geometric genus $7$ has two different plane models of degree $6$. The second one is obtained via a Cremona transformation through the $3$ non collinear double points. In our case however, the three singular points lie on a line $L=\{\ell=0\}$,  and we have a distinguished model. By adjunction theory, the rational map 
 $$C' \dasharrow C \subset \PP^6$$
 to the canonical model of such a sextic is defined by the  7-dimensional space of cubics through the double points. Since the double points lie on a line,  this space is spanned by  6 quadrics multiplied by $\ell$ and a further cubic, which intersects  the line in the three points. We will use all of this information in our construction.

\section{Construction procedure} 

Consider Euler parameters on $SO(3)$, i.e.,  the parametrization by unit quaternions, or, more conveniently, the identification of $\PP^3 \cong\overline{SO(3)} \subset \PP^9$ via the 2-uple embedding
$$\PP^3 \hookrightarrow\PP^9, \quad (s_0,s_1,s_2,s_3) \mapsto (A(s),t) \in \PP^9$$
defined  by 
$$A(s)=\
\begin{pmatrix}
s_{0}^{2}-s_{1}^{2}-s_{2}^{2}+s_{3}^{2}&
      -s_{0} s_{1}+s_{2} s_{3}&
      -s_{0} s_{2}-s_{1} s_{3}\\
      -s_{0} s_{1}-s_{2} s_{3}&
      -s_{0}^{2}+s_{1}^{2}-s_{2}^{2}+s_{3}^{2}&
      -s_{1} s_{2}+s_{0} s_{3}\\
      -s_{0} s_{2}+s_{1} s_{3}&
      -s_{1} s_{2}-s_{0} s_{3}&
      -s_{0}^{2}-s_{1}^{2}+s_{2}^{2}+s_{3}^{2}\\
      \end{pmatrix}$$
and $t=s_0^2+s_1^2+s_2^2+s_3^2$. 

In particular, we note that every curve in 
$\overline{SO(3)} \subset  \PP^9$ has even degree. 
This explains why even degree  motion curves are more likely than  odd degree curves. Odd degree can only occur if $C$ intersects the 2-dimensional vertex within the 6-dimensional  cone $\overline{SE(3)}$.\medskip

Utilizing the special geometry of the genus 7 curves in $\overline{SE(3)}$ as a hint, we will construct exceptional mechanisms in three steps. We say that a curve $C  \subset SE(3)$ supports the leg defined by $p^{(i)},q^{(i)}$ if the equation $eq^{(i)}$ is satisfied for all points $(A,b) \in C$.

\begin{enumerate}
\item Choose a plane $\PP^2 \subset \PP^3$ through the point $(0:0:0:1)$ corresponding to the identity in $SO(3)$ and choose three of the desired six legs $p^{(1)},q^{(1)},\ldots, p^{(3)},q^{(3)}$ freely.
\item Compute the family of plane sextics $C' \subset \PP^2 \subset \PP^3 \cong \overline{SO(3)}$ together with an extension to a birational map $C' \to C \subset \overline{ SE(3)} \subset \PP^{12}$ such that the equations $eq^{(1)},eq^{(2)}$ and $eq^{(3)}$ are satisfied for all points on $C$. It turns out that this family depends unirationally on three parameters.
\item Choose the parameters such that $C$  supports further three legs. 
\end{enumerate}

To start, we can write the equation
$eq^{(i)}$ homogeneously in terms of Euler parameters as
$$2t^2\langle p^{(i)},q^{(i)}\rangle+2\langle Ap^{(i)},tq^{(i)}\rangle-2\langle Ap^{(i)},b\rangle+2\langle b,tq^{(i)}\rangle- \parallel b\parallel^2=0.$$
Since $C \subset \overline{SE(3)}\cap \PP^6$ should be the canonical embedding,  the entries of
$A$, $t$ and the $b_i$ are given by cubic forms on our $\PP^2$ where the $a_{ij}$ and $t$ have the common factor $\ell$, and the equations evaluate to an equation of the plane sextic. The difference $eq^{(i)}-eq^{(j)}$ has $\ell$ as common factor. Division by $\ell$ leaves
us with a system of equations of degree five on the plane. We may regard this as a linear
equation 
$$\begin{pmatrix}
\alpha_0& \alpha_1&\alpha_2& \beta\cr \end{pmatrix} 
\begin{pmatrix}b_0\cr{b_1}\cr{b_2}\cr{\ell}\end{pmatrix}=0$$
for $b_0,b_1,b_2$ and $\ell$  where the coefficients $\alpha_i$ and $\beta$ are quadrics and a quartic which depend on the choice of the two legs. Taking two differences, we obtain a system of equations
$$\begin{pmatrix}
\alpha_0& \alpha_1&\alpha_2& \beta\cr 
\alpha_0'& \alpha_1'&\alpha_2'& \beta'\cr 
\end{pmatrix} 
\begin{pmatrix}b_0\cr{b_1}\cr{b_2}\cr{\ell}\end{pmatrix}=0$$
and any solution substituted in one of the original equations will give a plane sextic which supports the first three legs.

\bigskip
\noindent{\bf Step 1.} 
We will choose the first three legs and a plane $V(s)=\PP^2 \subset \PP^3$ through the point 
$(0:0:0:1)$ corresponding to $id \in SO(3)$ with integral coordinates.
Then $S=\QQ[s_1,s_2,s_3]\cong \QQ[s_0,\ldots,s_3]/\langle s \rangle$, the coordinate ring of $\PP^2$, is defined over $\QQ$. 

\bigskip
\noindent{\bf Step 2.}
The solutions $(b_0,b_1,b_2,\ell) \in S^3\oplus S(-2)$ are elements of degree $3$ in the kernel of the homomorphism
$ S^3\oplus S(-2) \to S^2(2)  $
defined by
 $$\begin{pmatrix}
\alpha_0& \alpha_1&\alpha_2& \beta\cr 
\alpha_0'& \alpha_1'&\alpha_2'& \beta'\cr 
\end{pmatrix}.$$

For general choices (at least experimentally),  the kernel is the image of a map
$$ S(-2) \oplus S(-3) \oplus S(-4) \to S^3\oplus S(-2)$$
 defined by a syzygy matrix
$$ 
\begin{pmatrix} 
a_0&b_0'& c_0\cr 
a_1&b_1'& c_1\cr
a_2&b_2'& c_2\cr
0 &\ell& a_4\cr
\end{pmatrix},
$$
whose entries are homogeneous forms $a_i,b_i', c_i$ of degree 2, 3 and 4 and a linear form $\ell$. In particular, $L=V(\ell) \subset \PP^2$ is uniquely determined by the choices in step 1.
For $b$ we have a three parameter family of solutions
$b_i=b_i' + ga_i$
where $g=xs_1+ys_2+zs_3$ is a linear form with three yet unknown coefficients $x,y,z$.
The family of sextics is now given by 
$$f=2t^2\langle p^{(1)},q^{(1)}\rangle+2\langle Ap^{(1)},tq^{(1)}\rangle-2\langle Ap^{(1)},b\rangle+2\langle b,tq^{(1)}\rangle- b_0^2-b_1^2-b_2^2$$
which is  a homogeneous polynomial in $s_1,s_2,s_3$ with coefficients in $\QQ[x,y,z]$.

\bigskip
\noindent{\bf Step 3a.} We now turn things  around and ask for the condition that the curve $C\subset \overline{SE(3)}$ defined by $f$ supports three further legs. Equations for a leg $(p,q)$ are obtained by asking that
$f$ and $f_{pq} \in \QQ[p,q,x,y,z][s_1,s_2,s_3]$ with
$$f_{pq}=2t^2\langle p,q\rangle+2\langle Ap,tq\rangle-2\langle Ap,b\rangle+2\langle b,tq\rangle- \parallel b\parallel^2$$
are proportional as homogeneous polynomials in $s_1,s_2,s_3$.

This yields an ideal $I \subset  \QQ[p,q,x,y,z]=\QQ[p_0,p_1,p_2,q_0,q_1,q_2,x,y,z]$  whose zero set
is the collection of legs. Saturating with the ideal of the known legs, gives a simplified ideal
 $J \subset \QQ[p,q,x,y,z]$. A Gr\"obner basis computation shows that the zero loci $X=V(J)\subset \AAA^9$ is a surface of degree 18. We are interested
 in values $\bar x,\bar y,\bar z$ such that
 $$X \cap V(x-\bar x,y-\bar y,z-\bar z)$$ consists of three points. In other words, we are looking for a triple point of the projection 
 $$ X \dasharrow Z \subset \AAA^3, \;(p,q,x,y,z) \mapsto (x,y,z).$$
 Since $X$ is a surface, we expect a finite number of triple points. 

\bigskip 
\noindent{\bf Step 3b.} The ideal $J$ is generated by linear and quadratic polynomials in $p,q$ with coefficients in $\QQ[x,y,z]$. 
We take these polynomials and all multiples with any of the variables $p,q$.
The coefficient matrix with respect to monomials in $p,q$ up to degree 3 gives us a matrix with entries
in $\QQ[x,y,z]$. This matrix  drops rank by 3 at the triple point loci of the projection, we utilize this   to compute the triple point loci.

\section{Existence over a number field}\label{number field}

For three legs with integral coordinates of moderate height, the computation of the ideal $J$ can be done over the integers. 
But to compute the triple points, that is the ideal of minors in Step 3b,  is out of reach over $\QQ$. Instead  we pass first to various finite fields $\FF_p$ where this computation is possible. 
The degree of the triple point locus is 210 for nearly all primes. \medskip

 It turns out that among the coefficients $x,y,z$ the coefficient $z$ plays a particular role: Scanning through the number of triple points on the planes defined over $\FF_p$, that is $z=\overline z$ for $\overline z \in \FF_p$, we find for most primes four distinguished values $\overline z\in  \FF_p$ which contain a lot of triple points. Three planes contain 28 triple points with coordinates in $\overline \FF_p$, the fourth 14 triple points. There are 112 further triple points.
 
It turns out that $X$ is always reducible: Its projective closure $\overline X$ is
$$\overline X= \PP^2 \cup \PP^2 \cup \PP^2 \cup Y \subset \PP^9.$$
The residual surface $Y$ is a singular
 conic bundle over $\PP^1$. One can check this computationally using the special shape of the syzygies of $Y$. We do not explain this  computation here.
Each of the image of a $\PP^2$ contains 28 triple points, while the plane corresponding to the 14 triple points arises as follows:

Two of the conics get identified under the projection, and the projection of $Y$ has 14 triple points along the image conic. We decided to go for a point on this conic. The other family of triple points will most likely lead to number fields of degree 112, and a further extension of degree 6 for the coordinates of the three preimage points. We do not study this second  family any further here.

Using the LLL-algorithm \cite{LLL} or Wang's rational conversion \cite{Wang}, we can find the equation in $\QQ[x,y,z] $ for the special plane, for  the conic within this plane, and finally for the septic, which cuts the conic in the 14 triple points.
Thus we obtain number fields $\QQ \subset K \subset L$ of degrees $[K:\QQ] \le 14$ and $[L:K] \le 2$,  where $K$ contains the coordinates 
$x,y$ and the possibly quadratic extension $L \supset K$ contains all coordinates of the legs. The degree $[L:K]\le 2$, because the preimage point of the triple point, which does not lie on the union of the two conics is distinguished. The conics themselves  are just conjugated for most choices of three initial legs.

The specific example below is obtained by choosing the legs $[p^{(i)} ,q^{(i)}]$ for $i=1,2,3$ as
$$\Bigg[
\begin{pmatrix}  0 \cr 0 \cr 0 \end{pmatrix},\begin{pmatrix}  3 \cr 0 \cr 0 \end{pmatrix}\Bigg] ,\;\Bigg[ \begin{pmatrix} -2 \cr -2 \cr -1 \end{pmatrix},\begin{pmatrix} -3 \cr 3 \cr -3 \end{pmatrix}\Bigg]  ,\; \Bigg[\begin{pmatrix} 2 \cr 1 \cr -2 \end{pmatrix},\begin{pmatrix} 3 \cr 1 \cr 1 \end{pmatrix} \Bigg],$$
and by taking $\PP^2=V(s_0-s_1) \subset \PP^3$ for the plane.
We  have choosen this example out of many others since the plot of the plane sextic, looked nice to us, see Figure 4 below.

From our computations we obtain the equations for the coefficients $$z=\frac{90}{47},$$ the conic

$$Q=x^2+\frac{214}{181}xy-\frac{377}{181}y^2+\frac{266242}{25521}x
-\frac{499174}{25521}y-\frac{45419195}{1199487}$$
and the septic, which we display below. 

{\small
\begin{eqnarray*}
  S=&   40529281858604847902260745656467021 x y^{6}\\
&  -13003502555830472808716754914961492 y^{7}\\
&   +1823658874752524908287468326069696796 xy^{5}\\
&   -1345103945703671198677310366396663022 y^{6}\\
&   +21057817082494339849424039333426043801 x y^{4}\\
&   -21097516746760090326807604348970652732 y^{5}\\
&   +104146801533327515594541432822537162912 x y^{3}\\ 
&   -147061992679082503452275851430863761654 y^{4}\\
&   +261417542057521780350077326987208938863 x y^{2}\\
&   -573573843341002730574525811534222498236 y^{3}\\
&   +137627172885093280917715345045703910468 x y\\
&   -1174127497815997633825349561153169116826 y^{2}\\
&   -1021528567727081989434308426376376661 x\\
&   -1580629500417030142273626496530002994868 y\\
&   -1752162935565295197876926071913693437570 \; .
\end{eqnarray*}
}

The septic equation describes $x$ as a rational function of degree 7 of $y$ and substituting this into the quadric, yields the polynomial which defines the field extension $K$.

\section{Real solutions} To prove the existence of a non-empty open set of initial choices as in Step~1 such that the number fields $K \subset L$ have a real embedding, it suffices to establish this for a single example. Figure 3 
shows a plot of  the quadric and the septic of our specific example from Section \ref{number field}.
There are  6 real intersection points, out of which five are visible in the plot. We choose
$$(x,y)\approx (2.2072318327235894538, -7.4995700484259023117).$$
The corresponding plane sextic has indeed three double points, all real with real tangents, as we can see in the plot of the curve in an affine chart given in Figure 4. 

\begin{figure}[h]\label{QS}
\begin{center}
\includegraphics{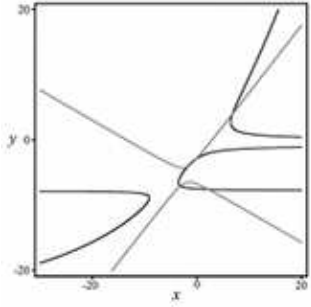}
\end{center}
\caption{A quadric and a septic with 6 real intersection points (5 are visible).}
\end{figure}

\begin{figure}[here]\label{sextic}
\begin{center}\includegraphics{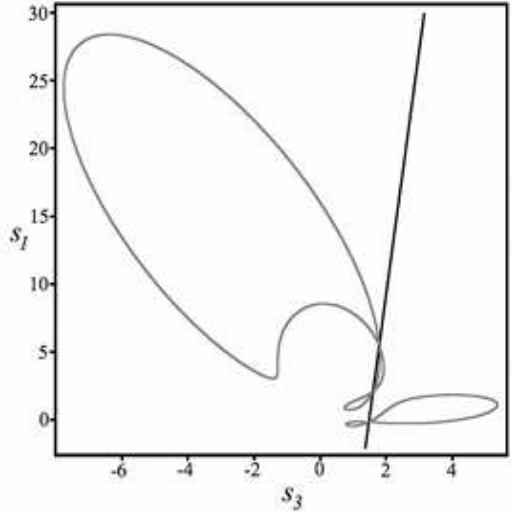}\end{center}
\caption{Plane sextic with three double points on a line}
\end{figure}

\noindent
We conclude from the picture that this specific  mechanism will have two disjoint connected components over $\RR$. 

Next, we compute the remaining legs $[p^{(i)},q^{(i)}]$ for $i=4,5,6$ using the ideal $J$ from Step 3a. They turn out to be real as well. The approximate values are
$$\Bigg[\begin{pmatrix}3.16918195183411\cr  2.10535216460007\cr  0.09679237835995\end{pmatrix},
\begin{pmatrix} 4.06043490691685 \cr 2.08171150266153\cr  3.29773800483749 \end{pmatrix}\Bigg]$$
$$\Bigg[\begin{pmatrix} 2.21062671304078,\cr 1.14679692580674 \cr -2.1941220884757 \end{pmatrix}, 
\begin{pmatrix} 3.10187966812352\cr  1.1231562638682 \cr  1.00682353800184\end{pmatrix} \Bigg ]
$$
  and 
$$\Bigg[\begin{pmatrix}3.5973027547094 \cr  3.5973027547094 \cr   1.79865137735471\end{pmatrix} , \begin{pmatrix}5.94796386794961\cr  4.40897636315916 \cr 1.4739819339748\end{pmatrix} \Bigg].$$

 \begin{figure}\label{positions}
\begin{tabular}{ccc}   
\includegraphics{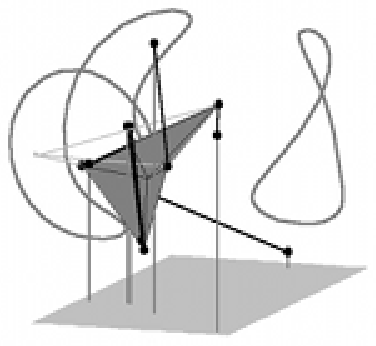} &
\includegraphics{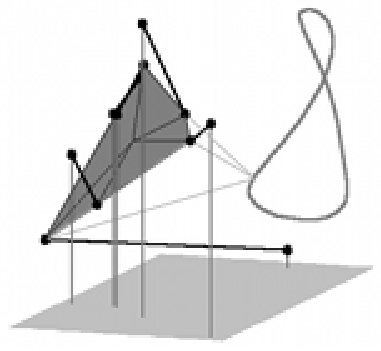} & 
\includegraphics{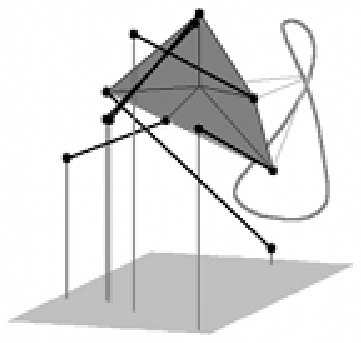} \\ 
\includegraphics{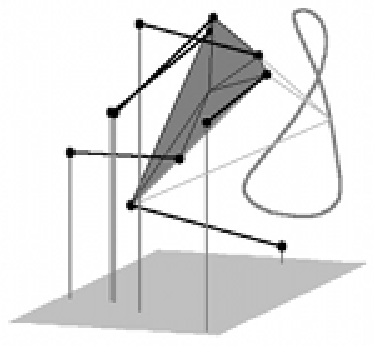} &
\includegraphics{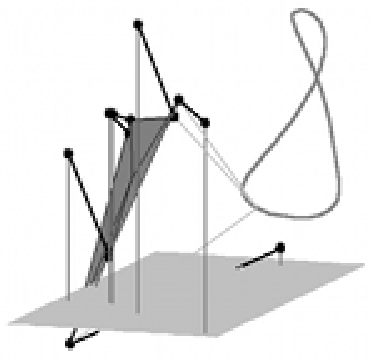} &
\includegraphics{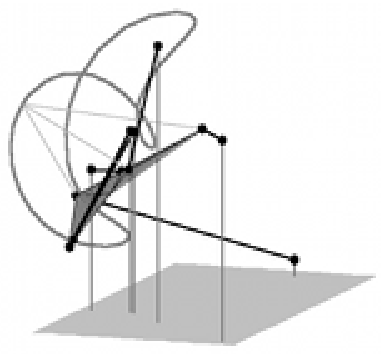} \\
\includegraphics{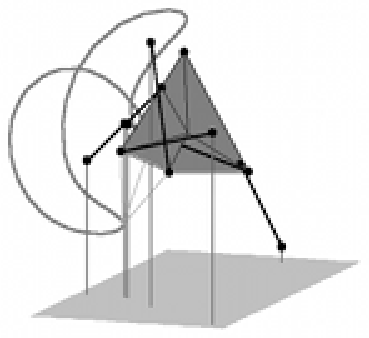} &
\includegraphics{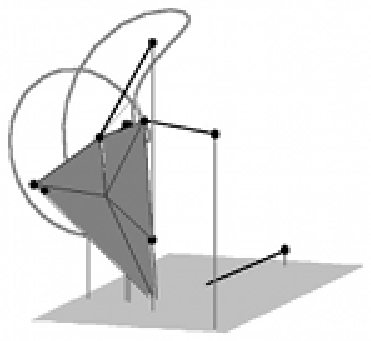} &
\includegraphics{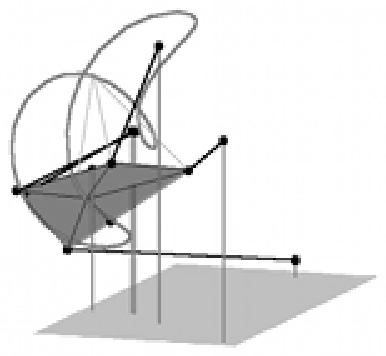} \\  
\end{tabular}
\caption{An exceptional Stewart-Gough platform. The
upper left is the rigid identity  position in the assembly space. The remaining pictures show the mechanism in different positions on both connected  components of the motion curve.
}
\end{figure}

 \noindent
 However, in this position the mechanism is rigid, because the point corresponding to $id \in SE(3)$ does not lie on the curve of motion. To obtain a start positions on each connected component, we apply the motion given by one element on each  connected component. Of course, the $q^{(i)}$ do not change. The coordinates of the $p^{(i)}$ for a start position
on left hand side component of in Figure \ref{positions} are  given by the 6-tuple of points with aproximate coordinates:
$$ p^{(1)} = \begin{pmatrix}  0.155967547996383 \\ -0.330589644955676 \\ -0.915719119302685 \end{pmatrix}, \ \ \
   p^{(2)} = \begin{pmatrix} -1.89049475348576 \\ -2.26381911010301 \\ -1.95235593308989 \end{pmatrix}, $$
$$ p^{(3)} = \begin{pmatrix}  3.02807585569042 \\ -1.02097882111392 \\ -1.43950977817185 \end{pmatrix}, \ \ \
   p^{(4)} = \begin{pmatrix}  3.54839453008461 \\  0.988949245364636  \\ 0.195967470045645  \end{pmatrix}, $$
$$ p^{(5)} = \begin{pmatrix}  3.34246862886166 \\  -1.05120251881891 \\  -1.50147789095895 \end{pmatrix}, \ \ \ 
   p^{(6)} = \begin{pmatrix}  3.83683978692307 \\ 3.14661620544245 \\ 0.948829119787298 \end{pmatrix}. $$
The right hand side component has a start position:
$$p^{(1)} = \begin{pmatrix} 0.651863003594838 \\ -0.291286328884768 \\ -1.84428994846972 \end{pmatrix},\ \ \ 
p^{(2)} = \begin{pmatrix}  -1.17539378163635 \\  -2.15417982581569 \\ -3.32441172045981 \end{pmatrix}, $$
$$p^{(3)} = \begin{pmatrix}-0.0897065626926926 \\ 0.326427753808606 \\ 0.99622125383145 \end{pmatrix}, \ \ \
p^{(4)} = \begin{pmatrix} 1.82686560387732 \\ 2.13947602458764  \\ 0.838312156777211 \end{pmatrix}, $$
$$p^{(5)} = \begin{pmatrix}-0.125381784762398 \\ 0.392226097496023  \\ 1.30926149281988 \end{pmatrix},\ \ \
p^{(6)} = \begin{pmatrix} 3.93846094695261 \\ 3.05940963208061 \\ 0.817933114521703 \end{pmatrix}. $$

Using a singular value decomposition version of path tracking due to Charles Wampler, one can trace out the motion curve numerically. 

\section{Open questions}
The family of mechanisms constructed depends on the $3\cdot 6$ parameters of the first three legs,
and $2$ parameters for the $\PP^2 \subset \PP^3 \cong \overline{ SO(3)}$ which passes through the point
$id \in SO(3)$. Hence altogether we  have 20 parameters, and  the corresponding family $M \subset \AAA^{36}$
has codimension 16. Considering equivalence classes up to translation, rotation, and dilation subtracts $7$ dimensions, thus we have a family of dimension $13$. Coming back to the experimental exploration, we did not expect to find points
on component of codimension much larger than $12$.  \medskip

{\bf Problem 1.} Explain why we were so lucky to find examples over $\FF_3$ nevertheless. \medskip
 
 \noindent
One possible explanation is that there are indeed many different components of degree 12 mechanisms in $\AAA^{36}$. Each could contribute with the same amount to the number of $\FF_p$-rational mechanisms. 
For example, the symmetry among the six legs is broken. The two legs lying on the conics, seem to have always the same length. Already in case of two distinguished legs, we obtain ${6 \choose 2}= 15$ different components in the space of 6 ordered legs. However, this is still not enough to account for the difference to the data from the experimental exploration. \medskip

\noindent
The generators of the syzygies of the matrix  $$\begin{pmatrix}
\alpha_0& \alpha_1&\alpha_2& \beta\cr 
\alpha_0'& \alpha_1'&\alpha_2'& \beta'\cr 
\end{pmatrix} $$ in Step 2 are of lower degree than the syzygies of a general map $$ S^3\oplus S(-2) \to S^2(2).  $$ 

{\bf Problem 2.} What is the reason for the low degree syzygies? \medskip

\noindent
The surface $\overline X= \PP^2 \cup \PP^2 \cup \PP^2 \cup Y$ of extra legs has a very special structure. \medskip

{\bf Problem 3.} Explain, $\dim X=2$ and the special structure of the surface $X$. In particular,  explain the conic bundle structure of $Y$. \medskip

{\bf Problem 4.} Why are two conics of $Y$  identified under the projection $X \dasharrow Z \subset \AAA^3$? \medskip

{\bf Problem 5.} The legs with coordinates on the two conics have the same length. Why? \medskip 

\noindent
It is not clear whether or not the resulting family of mechanisms is unirational.  \medskip 

{\bf Problem 6.} Determine whether this family of mechanism is unirational.  \medskip 

\noindent
For a non-unirational irreducible variety $M$ defined over $\QQ$ or $\RR$ the set 
 of real points $M(\RR)$ can have several connected components. 

{\bf Problem 7.} If the family of mechanisms is not unirational, determine the number of its connected components over the reals.\medskip

\noindent
Since we do not have a unirational parametrization of the family, it is not easy to change the mechanisms continuously. Perhaps, it is possible to set-up a  path tracking approach following the ideas of Sommese and Wampler \cite{SW}. 

{\bf Problem 8.} Vary the parameters of the mechanisms using homotopy methods! \medskip

{\bf Problem 9.} Compute the other 112 triple points and a corresponding mechanism! \medskip

{\bf Problem 10.} Study other families of mechanisms!

\end{document}